\documentclass[11pt,oneside]{amsart}
\usepackage {amsmath, amssymb,   epsfig, hyperref , enumerate,  color, graphicx}
\usepackage[utf8]{inputenc}
\usepackage[english]{babel}
\usepackage[paper=a4paper]{geometry}
  
\numberwithin{equation}{section}

\newcommand{\C}{\mathbb{C}}

\renewcommand{\P}{\mathbb{P}}
\newcommand{\e}{\varepsilon}

\newcommand{\fr}{\partial}
 \newcommand{\set}[1]{\left\{#1\right\}}
\newcommand{\norm}[1]{\left\Vert#1\right\Vert}
\newcommand{\abs}[1]{\left\vert#1\right\vert}
\newcommand{\cd}{{\mathbb{C}^2}}
\newcommand{\pd}{{\mathbb{P}^2}}
\newcommand{\pu}{{\mathbb{P}^1}}
\newcommand{\pk}{{\mathbb{P}^k}}

\newcommand{\rest}[1]{ \arrowvert_{#1}}

\newcommand{\unsur}[1]{\frac{1}{#1}}

\newcommand{\lrpar}[1]{\left(#1\right)}

\newcommand{\inv}{^{-1}}

\newcommand{\la}{\lambda}

\newcommand{\tbif}{T_{\mathrm{bif}}}
\newcommand{\mubif}{\mu_{\mathrm{bif}}}

\newcommand{\poly}{\mathcal{P}}
\newcommand{\modu}{\mathcal{M}}

\newcommand{\hot}{\mathrm{h.o.t.}}
\newcommand{\loc}{\mathrm{loc}}

\DeclareMathOperator{\supp}{Supp}
\DeclareMathOperator{\vol}{Vol}
\DeclareMathOperator{\Int}{Int}

\theoremstyle{theorem}

\newtheorem{prop}{Proposition} [section]
\newtheorem{thm}[prop] {Theorem}

\newtheorem{cor}[prop]{Corollary}

\newtheorem{question}[prop]{Question}

\theoremstyle{definition}

\begin{document}


\title{Geometric methods in holomorphic dynamics}

\author{Romain Dujardin}
\address{Sorbonne Universit\'e, Laboratoire de Probabilit\'es, Statistique et Mod\'elisation,    4 place Jussieu, 75005 Paris, France}
\email{romain.dujardin@sorbonne-universite.fr}

\begin{abstract}
In this note we review a selection of   contemporary  research themes  in  
 holomorphic dynamics. 
  The main topics that will be discussed are: geometric (laminar and woven) currents and their applications, 
bifurcation theory  in one and several variables,  
and the problem of wandering  Fatou components.
\end{abstract}

\maketitle

Holomorphic dynamics   was once part of classical complex analysis, 
but since its rebirth in the  1980's it   keeps enlarging its scope, integrating 
new ideas and     developing  new interactions. 
 Some main tendencies of contemporary holomorphic dynamics are 
 the convergence between its one and higher dimensional aspects and 
 its ever deeper interconnection with algebraic and  arithmetic dynamics.  
As a consequence, there is an endless diversification of  the available mathematical techniques.
 Besides the classical methods  from 
 dynamics and   complex analysis, its  modern  
  toolbox now  comprises sophisticated  tools and ideas imported from:  complex geometry, 
 pluripotential theory (and its latest advances for currents of higher bidegree), 
  algebraic geometry and commutative algebra, 
  non-Archimedean analysis and geometry,
  arithmetic geometry (in particular arithmetic equidistribution theory),
  Teichmüller theory,  geometric group theory,  etc. 
Conversely, each of these domains benefits from its interaction with holomorphic dynamics, by gaining new problems and examples. 
Many  (though not all!) of these connections were reported in recent 
ICM's \cite{buff-cheritat, lyubich_icm, oguiso_icm, cantat_icm, demarco_icm, dinh_icm}. 
Our purpose here is to  present  a few contemporary research themes whose common thread --if one were to find one-- is an    
emphasis on ``soft''  geometric techniques, such as  the basic    geometry of analytic subsets in $\C^n$.
These represent only a tiny piece of the domain, reflecting of course the author's own taste and research interests. 
The main topics that will be discussed are: geometric currents, 
bifurcation theory, and the problem of wandering  Fatou components. The reader will soon notice
 that these three subjects  are  largely   interrelated.  Many  open problems have also 
  been included, as a motivation for future 
 investigations.

Let us describe in more detail the contents of this paper. Section~\ref{sec:currents} 
is a short  survey on positive closed 
currents with  ``geometric structure''. The use of geometric currents 
 in holomorphic dynamics was pioneered by Bedford, Lyubich and Smillie in their 
 seminal work~\cite{bls} on complex Hénon maps. Since then they have turned  into a 
 very versatile tool, with many applications. Here we intend to  give the   
  flavor of a few specific  results 
and how they are used in dynamical problems, so this part of the paper 
 will be a bit more     technical than the remaining sections.

Holomorphic dynamics is equally about the dynamics  of a holomorphic map $f$ and 
about the evolution of  this dynamical 
behavior   when $f$ depends on certain parameters. 
 The basic  stability/bifurcation theory of rational maps 
 in one variable was designed by  Mañé, Sad, Sullivan and Lyubich \cite{mss, lyubich_stability1, lyubich_stability2} in the 1980's, 
who showed 
that one-dimensional rational maps are generically structurally stable, using surprisingly elementary arguments. For the  
  quadratic family $z^2+c$, $c\in \C$, the bifurcation locus is the celebrated Mandelbrot set, whose intricate structure was 
  thoroughly studied since then, using a variety of combinatorial and geometric methods. 
  This research area was profoundly renewed   in the 2000's
by the systematic investigation  of higher dimensional phenomena, and in particular with  the 
 introduction of bifurcation currents by DeMarco~\cite{demarco_current}.
The  bifurcation theory of holomorphic dynamical systems is nowadays a very active research domain, 
and a meeting point between the communities of one and several variable dynamicists. 
We relate this continuing story in  Section~\ref{sec:bifurcation}.
  
 Finally, one recent breakthrough is the construction of wandering Fatou components in higher dimensional polynomial dynamics, which at the same time solves an old problem and   raises many  questions. We review these recent developments in Section~\ref{sec:wandering}. 
 
Let us conclude this introduction with a little notice: some important  theorems will be mentioned only in passing, 
while other are isolated within  numbered environments:
this is meant to keep the reading flow, not to
 reflect a hierarchy of importance.  Likewise, the list of references is already quite  long, but not exhaustive, 
 and we apologize in advance for 
 any serious omission.  

\bigskip 
\noindent{\bf Acknowledgments.} It is a great pleasure to thank my collaborators  
 Eric Bedford, Pierre Berger, Serge Cantat,  Bertrand Deroin, Jeff Diller, 
 Charles Favre,  Vincent Guedj,  and Misha Lyubich 
for so many discussions about the mathematical  themes presented here, and more generally my colleagues 
from the holomorphic dynamics community for maintaining such  a friendly atmosphere over the years. 
 Special thanks to Charles Favre and Thomas Gauthier for their helpful comments on  this paper.

\section{Geometric currents}\label{sec:currents}~

\subsection{Definitions} 
This part assumes some familiarity with positive currents and pluripotential theory (see e.g. Demailly~\cite{demailly_agbook} for basics). 
All the definitions here are local so we
 work in some bounded open set $\Omega\subset \C^k$. Let $T$ be a positive
closed current of bidimension $(p,p)$ in $\Omega$. 
Following Bedford, Lyubich and Smillie~\cite{bls}, we  say that $T$ is \emph{locally uniformly laminar} if there 
exists a lamination by complex submanifolds of dimension $p$ embedded in $\Omega$ such that the restriction 
of $T$ to any flow box $B$ of the lamination is of the form 
\begin{equation}\label{eq:UL}
T\rest{B} = \int_\tau [\Delta_t] d\nu(t).
\end{equation} Here $\tau$ is a global  transversal in the flow box $B$, the $
\Delta_t$ are the plaques of the lamination in the flow box, and $\nu$ is a positive measure on $\tau$. 
The word ``uniform'' here refers to the local uniformity of the geometry of the plaques $\Delta_t$ 
We say that $T$ is \emph{laminar} if there exists a sequence of open subsets $\Omega_k$, 
together with a sequence of currents $T_k$, 
locally uniformly laminar in $\Omega_k$, such that $T_k$ increases to $T$. The $\Omega_k$ should be thought of as a union of many small polydisks, whose complement has a small mass. 
The key word in the definition  is ``increases'': 
intuitively, this definition should be understood as follows: $T_k$ represents all the disks contained in  $T$ of 
some given size (say $2^{-k}$); then, to $T_k$   we add $T_{k+1}- T_k$ which is made of       
 disks of size $2^{-(k+1)}$ (which may have non-empty boundary in $\Omega_k$, but form a lamination in 
 $\Omega_{k+1}\subset \Omega_k$), and so on. The sequence $T_k$ is not canonical, and has to be understood as the choice of a  ``representation'' of $T$ as a laminar current. 
From this we can deduce another    representation of $T$ as an 
 integral over an abstract family of   {compatible} holomorphic disks: 
\begin{equation}\label{eq:L}
T = \int_{\mathcal A} [D_\alpha] d\mu(\alpha). 
\end{equation} Here 
 \emph{compatible} means that two disks can only intersect along some relatively open subset, but there is no  further 
 restriction on the geometry of the $D_\alpha$. Even if this definition is rather 
 restrictive, it can lead to   pathological examples, and for dynamical applications we 
 will have to constrain it further (see  the notion of ``strongly approximable'' current below).
   
It was observed by Dinh~\cite{dinh_laminaire} that in many situations it is more natural to let the disks 
admit non-trivial  intersections. One then defines \emph{uniformly woven}  currents by
 replacing ``lamination'' by ``web'' in 
\eqref{eq:UL}, where a web   is locally given by a family of disks of dimension~$p$ 
with uniformly bounded volume, or more generally 
a family of holomorphic chains of dimension~$p$ 
with uniformly bounded volume (any such family is pre-compact for the Hausdorff topology, so it makes sense to define a measure on a set of such disks). Then,  \emph{woven}  currents  are defined 
from   uniformly woven ones as in the laminar  case. A difference between laminar and woven currents is that 
in the woven case the measures in~\eqref{eq:UL} and~\eqref{eq:L} are not determined by $T$ (e.g.  the 
 standard Kähler form in $\cd$ admits several representations as a uniformly woven current), 
 so a woven current 
 has to be thought of as ``marked'' by such a measure $\mu$. It is not completely obvious to show that not every 
 positive closed current is woven:   we leave this as an exercise to the reader! 
 
There is no unified reference for the basic properties of laminar and woven currents. Besides \cite{bls} and \cite{dinh_laminaire}, the  information in this paragraph was extracted from various papers, notably by De Thélin and the author~\cite{lamin, isect, fatou, ddg3, dethelin, dethelin2}.
 In the following we use the word \emph{geometric} as a synonym  of  ``laminar or woven".
  
\subsection{Construction and  approximation}\label{subs:approximation}

Positive closed currents   often appear as limits of  sequences of 
normalized currents of integration. 
Furthermore, by  a classical theorem of Lelong, any positive closed current of bidegree (1,1) 
is locally of this form. 
In this section we explain how under appropriate hypotheses, 
a geometric structure can   be extracted from such an approximation.  

Still working locally in some open set $\Omega\subset \C^k$, endowed with its standard Hermitian structure, we 
say that a submanifold $V$ of dimension $p$ in 
$\Omega$ \emph{has size $r$ at} $x\in V$, if it contains  a graph over a ball of radius $r$ of 
its tangent space  $T_xV$, relative to the orthogonal projection to $T_xV$, with slope (i.e. the norm of 
the  derivative of the graphing map) bounded by 1. In particular $V$ has no boundary in $B(x, cr)$ for some 
  constant $c$ depending only on $p$ and $k$. This notion of size makes sense in any compact
   complex manifold,  up to uniform constants, 
by choosing a finite covering by coordinate charts and  a Hermitian metric. 
Note that we may relax this definition by allowing $V$ to be an analytic set: then $V$ can have several irreducible components at $x$, some of which being of size $r$.

If $V$ is any submanifold (or subvariety) of $\Omega$, possibly with boundary, and $r>0$, 
we denote by $V^r$ the set of $x\in V$ such that 
$V$ has size $r$ at $x$. In this way we get a tautological decomposition: 
$V=V^r \cup (V\setminus V^r)$, which is reminiscent of the thin-thick decomposition of hyperbolic manifolds . 

Assume now that $V_n$ is a sequence of $p$-dimensional subvarieties of volume $v_n$, 
such that  $v_n\inv[V_n]$ converges to a positive closed current $T$.
 If $\vol(V_n^r)\geq v_n(1-\e(r))$ where 
$\e$ is a function independent of $n$  and   such that   $\e(r) \to 0$ as $r\to 0$,  
then one may extract a subsequence so that $v_n\inv[V_n^r]$ converges to a
    geometric current 
$T^r\leq T$ with the mass estimate 
$\mathbf{M}(T-T^r)\leq \e(r)$. This    
endows $T$    with  a geometric structure: if $p \leq k-2$
we obtain  a woven current and if $p = k-1$ this current is laminar. Indeed    
if $p = k-1$,  by the persistence of proper intersections,  
the limiting graphs cannot intersect non-trivially.  
(Note that when  $p\leq k-2$, intersections can appear at the limit even if the $V_n$ are 
submanifolds. Conversely, 
if in codimension 1 we allow the $V_n$ to admit 
 self-intersections, we    obtain woven currents also in this case.) 
 
A technically convenient option is to further assume 
 that the disks constituting $V_n^r$ are submanifolds (without boundary) in a 
 subdivision of $\Omega$ by cubes of size $cr$ 
 (for some   constant $c>0$). 
 This is consistent with the manner in which the $V_n^r$ are constructed in practice, and 
 the resulting definition is equivalent   (see \cite{fatou}). In this way the limiting currents 
  $T^r$ are uniformly geometric in the cubes of this subdivision.

 There are several easily checkable geometric and/or topological criteria ensuring this condition, which sometimes give an explicit bound on $\e(r)$:
\begin{itemize}
\item If $\psi:\C\to X$ is an entire curve in a projective manifold, then by Ahlfors' theory of covering surfaces, 
 for well chosen sequences $R_n\to\infty$, $V_n:=\psi(D(0, R_n))$ satisfies  $v_n\inv[\fr V_n] \to 0$ and 
 $\vol(V_n^r)\geq v_n(1-\e(r))$ for $\e(r) = O(r^2)$. Thus 
  the  cluster values of $v_n\inv[V_n]$ are closed woven currents; if in addition $\psi$ is injective and $\dim(X)=2$,  then 
 they are laminar (Bedford-Lyubich-Smillie~\cite{bls}, Cantat~\cite{cantat_K3}.
 \item If $V_n$ is a sequence of algebraic curves 
  in a projective surface whose geometric genus is $O(v_n)$, then 
 $\vol(V_n^r)\geq v_n(1-\e(r))$ for $\e(r) = O(r^2)$, therefore the limiting currents of $v_n\inv[V_n]$ are woven; under a mild additional condition on the singularities of $V_n$, they are laminar (Dujardin~\cite{lamin}).
\item If $\iota_n:\P^p\to X$ is a sequence of holomorphic mappings of generic degree 1 to a projective manifold $X$ of dimension $k> p$,  and $V_n   = \iota_n(\P^p)$   then 
 the limiting currents of $v_n\inv[V_n]$ are woven (Dinh~\cite{dinh_laminaire}). 
 In addition  $\e(r) = O(r^2)$ \cite{fatou}. 
 \item If $V_n$ is a sequence of smooth curves in the unit ball in $\cd$, whose genus is $O(v_n)$, then 
 the limiting currents of $v_n\inv[V_n]$ are laminar (De Thélin~\cite{dethelin}). A version of this result  
  in arbitrary dimension  is given by De Thélin in \cite{dethelin2}. 
\end{itemize}
In all these papers, the geometric structure is obtained by projecting $V_n$ in 
several directions and keeping  only from $V_n$ the 
 graphs over these directions with bounded diameter or volume. 
The bound $\e(r) = O(r^2)$   plays an important role in  applications as we shall see below.

 \subsection{Geometric intersection}
 
 The main interest of geometric currents is the possibility of a geometric interpretation of their wedge products. 
 This technique was introduced in \cite{bls},   and it was systematized and 
 generalized in several subsequent works. 
Such results are  so far  essentially available in  dimension 2; again
 since the problem is local  we work in some open set $\Omega\subset \cd$, say a ball. If $T_1$ and $T_2$ are closed positive (1,1) currents in $\Omega$, we say that the wedge product $T_1\wedge T_2$ is well defined  if $u_1\in L^1_\loc(T_2)$, where $u_i$ is a local potential of $T_i$, in which case we set 
 $T_1\wedge T_2 = dd^c(u_1T_2)$. This condition and the resulting wedge product 
 is actually symmetric in $T_1$ and $T_2$. We also say that  such a current is \emph{diffuse}
  if it gives no mass to curves. 
 
 For uniformly laminar and woven currents, geometric intersection is easy and   
basically follows from Fubini's theorem. Indeed, assume that $T_1$ and $T_2$ 
 are uniformly geometric (1,1)-currents in $\Omega$, which locally in 
 $\Omega$ admit the representation
 $T_i = \int[\Delta^i_t]d\nu_i(t)$. Then, if the wedge product $T_1\wedge T_2$ is well defined, 
  then locally we have that 
 \begin{equation}\label{eq:geometric}T_1\wedge T_2 = \int [\Delta^1_t\cap \Delta^2_s]d\nu_1(t)d\nu_2(s), 
\end{equation} where 
 $[\Delta^1_t\cap \Delta^2_s]$ is the sum of point masses at isolated  intersection points, counting 
 multiplicities (see \cite{isect, ddg2}).  In addition, if $T_1$ and $T_2$ are laminar and diffuse, non-transverse intersections 
 do not contribute to the integral so we can restrict to 
  transverse intersections. Note the intermediate ``semi-geometric intersection'' result 
\begin{equation}\label{eq:semi-geometric}
T_1\wedge T_2 = \int ([\Delta^1_t]\wedge T_2) d\nu_1(t), 
\end{equation} 
which makes sense for an arbitrary positive closed current $T_2$. 
 
Now assume that $T$ is a geometric positive closed 
current in $\Omega\subset \cd$  and $S$ is an arbitrary positive closed current in $\Omega$, such that the wedge product $S\wedge T$ is well-defined. 
We say that $T\wedge S$ is \emph{semi-geometric} if there is  a representation $T = \lim_{r\to 0} T^r$ 
as an increasing limit of uniformly geometric currents, such that 
  $T^r\wedge S$ increases to $T\wedge S$ as $r\to 0$. Thanks to~\eqref{eq:semi-geometric},  
 $T^r\wedge S$   admits a geometric interpretation.   
  If now $S$ itself is a geometric current, we say that the wedge product 
$T\wedge S$ is \emph{geometric} if there are representations $T^r\nearrow T$ and $S^r\nearrow S$
such that $T^r\wedge S^r$
(which has a geometric interpretation by~\eqref{eq:geometric})  
increases to $T\wedge S$.  

We say that a geometric current is \emph{strongly approximable} if there is a representation 
$T^r\nearrow T$ where $T^r$ is uniformly geometric  in a  subdivision $\Omega^r$ of $\Omega$ 
into cubes of size $r$, and $\e(r) = \mathbf{M}(T-T^r) = O(r^2)$. 
As we have seen in \S\ref{subs:approximation}, this estimate is commonly  
satisfied in practice. (Technically, some freedom 
on the choice of $\Omega^r$ is also necessary, but we do not dwell on this point.) 
The sharpest version of the geometric intersection theorem for geometric currents in dimension 2 is the following: 

\begin{thm}[Dujardin \cite{isect, birat, ddg2}] \label{thm:isect}
Let $S$ and $T$ be closed positive (1,1)  currents in $\Omega\subset \cd$, such that the wedge product 
$T\wedge S$ is well-defined. 
Assume that  $T$ is a strongly approximable  geometric current. Then, if $S$ has locally bounded potentials, or if
$T\wedge S$ gives no mass to pluripolar sets, then  $T\wedge S$ is semi-geometric. 
\end{thm}

A consequence of this theorem, which is often as useful as the result itself, is that if $T$   
was obtained 
as the limit of $v_n\inv[V_n]$   as in \S\ref{subs:approximation}, 
then $v_n\inv[V_n^r] \wedge S$ is close to 
$T\wedge S$ for small $r$ and  large $n$. 

Applying Theorem~\ref{thm:isect} to   $T\wedge S$ and $S\wedge T$, we get:

\begin{cor}
If in Theorem~\ref{thm:isect} both $S$ and $T$ are strongly approximable  geometric currents and 
$T\wedge S$ gives no mass to pluripolar sets, then  $T\wedge S$ is geometric.
\end{cor}

The main open problem at this stage is the extension of these results to higher dimensions. 

\begin{question}
Is there a version of Theorem~\ref{thm:isect} for geometric currents of arbitrary codimension?
\end{question}

  While  the case of uniformly geometric currents and the case where $T$ is of bidimension (1,1) follow
   without serious difficulties (see \cite{fatou} and \cite{higher}  for details),  the 
  general case remains a challenge so far.   The crucial mass estimate $\mathbf{M}(T-T^r) = O(r^2)$ is known to 
 hold  in some significant cases (see~\cite{fatou}), but it does not appear to be sufficient to conclude
for currents of arbitrary bidimension.
 
\subsection{Dynamical applications} \label{subs:lamin_applications}
The first application of laminar currents by Bedford, Lyubich and Smillie~\cite{bls} was to prove that certain intersections are non-empty. A typical example is the following: assume that  we are given 
an entire curve $\psi:\C\to X$  in some projective manifold, and let $T$ be a closed 
 current   obtained from $\psi$ by Ahlfors' construction.  Let $S$ be a current of 
bidegree (1,1) with bounded potentials. If we know that $\int T\wedge S>0$ (for instance, for cohomological reasons), then 
by Theorem~\ref{thm:isect}, this intersection is semi-geometric, therefore $S\rest{\psi(D(0, R_n))}$ is non-zero for large $n$. (A version of this result which does not appeal to laminarity was proved by Dinh and Sibony~\cite{dinh-sibony_automorphisms}.)
This fact (as well as some variants) plays an important role in the 
dynamics of automorphisms and birational maps on complex surfaces, where it is used as a tool to create intersections between 
stable and unstable manifolds. This is used  in~\cite{bls} to establish 
 that any saddle point belongs to the support of the maximal entropy 
measure; this technique  
 also appears in the work of Cantat, Favre, Lyubich, and the author~\cite{cantat_K3, tangencies, dmm, stiffness}, among others. Note also that the failure of 
Theorem~\ref{thm:isect} for unbounded potentials can be viewed as the main reason why  
the uniqueness of the measure of maximal entropy 
for general birational maps of surfaces remains an unsolved problem. 

Another use of geometric intersection, which was initiated in~\cite{birat},  concerns 
 the dynamical analysis of wedge products of 
dynamically defined currents. Indeed, suppose that $f$ is a self-map of some complex manifold $X$, 
and $f^n(L)$ is a sequence of 
iterated curves such that $d^{-n} f^n(L)$ converges to a geometric current $T$, with a control of the asymptotic geometry of $f^n(L)$ as in 
\S\ref{subs:approximation}. Assume also that $S$ is some invariant current of bidegree (1,1): $f^*S=dS$ and that $T\wedge S$ is a 
semi-geometric intersection. Then  for   large $n$,  the action  of   $  f^k $ on the bounded geometry part of 
  $d^{-n} [f^n(L)]\wedge S$  is a good approximation of the action of $f^k$ on $T\wedge S$, and 
    its expansion properties ``in the direction of $T$''
 can be analyzed geometrically  by ``soft'' methods, 
  such as counting disjoint disks of size $r$ and length-area estimates (see below Theorem~\ref{thm:dtgv} for 
  a worked out example). This idea was used in various contexts by De Thélin and others~\cite{birat, dethelin_selle, fatou, ddg3, dethelin-gauthier-vigny}.
  
  \subsection{Foliations}
Foliated Ahlfors  currents  
play an important role in the work of Brunella and McQuillan 
on singular holomorphic foliations (see e.g.~\cite{brunella}). 
Geometric intersection  has been 
 applied  in foliation theory to prove the vanishing of certain self-intersections. 
For a positive current  directed by a holomorphic foliation
 on a compact Kähler surface,   this vanishing  can in turn be used to infer dynamical   
 properties of the foliation such as 
the non-existence of invariant transverse measures (for closed currents)
 or the uniqueness of harmonic measures (for $dd^c$-closed currents),   according to a Hodge-theoretic 
 formalism for $dd^c$-closed currents devised by Forn\ae ss and Sibony~\cite{fornaess-sibony_laminations}. 
 Proving that the self-intersection of   harmonic currents directed by holomorphic foliations 
 vanishes  is a very difficult  problem in the 
 presence of singularities. On $\pd$  
 this can be treated by regularizing with global automorphisms,  
   the general case makes use of 
    the theory of densities of Dinh and Sibony 
  (see~\cite{dinh-nguyen-sibony}). Here we want to mention a more elementary-looking problem:  
  
 \begin{question}
 Does there exist a  diffuse (closed) uniformly laminar  current on $\pd$?
 \end{question}

 The expected answer to the question is  ``no'', since it is generally   expected that   
 there does not 
 exist  a Riemann surface  lamination embedded in $\pd$. The above question is supposed to be the ``easy 
 case'' of this deep conjecture (since it deals with laminations with transverse measures), and it admits a  
 straightforward approach: if $T$ is such a current, then $T\wedge T = 0$ because of the
  laminar structure, which is impossible on $\pd$. This approach works well as soon as  $T\wedge T$ is well-defined 
  in the sense of  pluripotential theory (but it doesn't for a curve!), or when the holonomy of the induced 
  lamination is Lipschitz~\cite{fornaess-sibony_laminations}. But in general the holonomy of a Riemann surface lamination in $\cd$ (that is, a holomorphic motion) is less regular and, surprisingly enough, 
   the problem is still open so far. (See Kaufmann~\cite{kaufmann} for a discussion of the higher dimensional case.)

%

\section{Bifurcation theory in one and several dimensions}\label{sec:bifurcation}~

Let $(f_\lambda)_{\lambda\in \Lambda}$ be a family
of rational maps on $\pu$ 
of degree $d$,   holomorphically parameterized by some complex manifold $\Lambda$. 
Then the  well-known Fatou-Julia decomposition of the phase space is mirrored by 
a stability-bifurcation dichotomy of the 
parameter space. The proper definition of  stability in this context was found 
simultaneously by Ma{\~n}é-Sad-Sullivan and Lyubich~\cite{mss, lyubich_stability1, lyubich_stability2}: 
the  family $(f_\lambda)_{\lambda\in \Lambda}$ 
is \emph{$J$-stable} over some domain $\Omega\subset \Lambda$
if one of the following equivalent conditions holds over $\Omega$:
\begin{enumerate}[(i)]\label{conditions}
\item the periodic points of  $(f_\la)$  do not collide, or equivalently the nature  
(attracting, repelling, indifferent) of each periodic point remains the same in the family; 
\item the Julia set $\lambda \mapsto J_\la$ moves continuously for the Hausdorff topology;
\item for any two parameters 
 $\la,\la'$ in $\Omega$, $f_\la\rest{J_\la}$  is topologically conjugate to $f_{\la'}\rest{J_{\la'}}$; 
 \item the orbits of the critical points $f_\lambda$ do not bifurcate. 
\end{enumerate}
The equivalence between these properties relies on  the notion of 
\emph{holomorphic motion} (also known as \emph{holomorphic families of injections})
of a subset of the Riemann sphere, and the 
the simple yet powerful idea of automatic  extension   of a holomorphic motion to its closure 
(the ``$\lambda$-lemma''). 
Condition (iv), together with the finiteness of the critical set, easily implies that in any such parameterized family 
$(f_\lambda)$, \emph{the stability locus is open and dense in $\Lambda$}. In other words, 
\emph{one-dimensional polynomial 
and rational  maps are generically stable}. 

For the emblematic   family  $f_\lambda(z) = z^2+\lambda$ of quadratic polynomials, the bifurcation locus 
is the boundary of the Mandelbrot set $M$ (connectivity locus). Even if its interior is empty, 
$\partial M$ is still quite 
large, as shown by the following famous result of Shishikura~\cite{shishikura}: \emph{$\partial M$ has Hausdorff 
dimension 2. } This property was extended to arbitrary families of rational maps by Tan Lei and McMullen 
\cite{tanlei, mcmullen}. The basic technical tool underlying Shishikura's theorem is the phenomenon of \emph{parabolic implosion}, which will also play an important role below. 
Note that  is still unknown whether $\fr M$ has zero or positive Lebesgue measure. 

This research area was   renewed in the last 20 years as the result of several tendencies: (1)~the use 
of positive closed currents, and (2)~the move  towards higher dimensions (both in dynamical and parameter 
spaces). In the next few pages we   review some of these developments; in particular we will see 
how these influential one-dimensional results translate to new settings. Lack of space prevents us from giving 
a complete treatment, and some important  results will barely be mentioned. Also, 
we do not discuss the profound connection  with arithmetic dynamics, for which the reader is referred e.g. to~\cite{demarco_icm}, nor bifurcations of Kleinian groups (see \cite{kleinbif, survey}).

\subsection{Bifurcation currents in one-dimensional dynamics}~
Let as above $(f_\lambda)_{\lambda\in \Lambda}$ be a holomorphic family of rational maps of degree $d$. 
The following addition to the list of equivalent conditions to stability was  found by  
 DeMarco~\cite{demarco_current2}:
\begin{enumerate} \label{condition5}
\item[(v)] the Lyapunov exponent of the unique measure of maximal entropy $\chi(\mu_{f_\lambda})$ is a pluriharmonic function of $\lambda$.  
\end{enumerate}
The \emph{bifurcation current} is then defined by $\tbif:= dd^c_\lambda \chi(\mu_{f_\lambda})$. 
 For the family of quadratic polynomials, $\tbif$ ($=\mubif$, see below) is the harmonic measure of the Mandelbrot set.  

The original definition of the bifurcation current in~\cite{demarco_current} can be interpreted geometrically as follows (see~\cite{prepercrit}).  Consider the fibered dynamical system in  $\Lambda\times \pu$
defined by $\hat f:(\lambda, z)
\mapsto (\lambda, f_\lambda(z))$. It admits a natural invariant current $\hat T$ of bidegree (1,1), 
satisfying $\hat f ^*\hat T  =d\hat T$, 
whose restriction  to a generic  vertical line  $\set{\lambda}\times \pu$ is the maximal entropy measure $\mu_{f_\lambda}$. 
Now, take a holomorphically moving (or ``marked'') point $\lambda \mapsto a(\lambda)$ in $\pu$, and denote by  
$\Gamma_a$ its graph in $\Lambda\times \pu$. If $\pi:\Lambda\times \pu\to\Lambda$ is 
the natural projection, we   obtain a current in $\Lambda$ associated to $a$ by slicing $\hat T$ by 
$\Gamma_{a}$ and projecting down to 
$\Lambda$: $T_a:= \pi_*(\hat T \wedge [\Gamma_{a}])$. 
If in a holomorphic family $(f_\lambda)$, the critical points are  marked by holomorphic 
functions $\lambda\mapsto c_i(\lambda)$ (this is always possible up to replacing $\Lambda$ by some branched cover), we thus
obtain the corresponding  bifurcation currents $T_{c_i}$. 
It turns out that  $\tbif  = \sum T_{c_i}$: 
 this follows from a variant of the Manning-Przytycki  formula for the Lyapunov exponent  
 $\chi(\mu_{f_\lambda})$, which in the case of polynomials writes 
 $$\chi(\mu_{f})  = \log d +\sum_i G_f(c_i),$$ where $G_f$ is the dynamical Green function (which satisfies 
 $dd^cG_f  = \mu_f$).  
  
Bifurcation currents have turned into a fundamental tool for exploring higher dimensional issues 
in parameter spaces. Here is a sample problem: 
consider a critically marked family $(f_\lambda, c_i(\lambda))$ and 
suppose that for some parameter $\lambda_0\in \Lambda$, the  critical point $c_1(\la)$ bifurcates at $
\lambda=\lambda_0$. Then a simple application of Montel's theorem shows that there is a sequence of parameters $
\lambda_n\to \lambda_0$ such that, for  $\lambda = \lambda_n$, $c_1$ is preperiodic. 
Now assume that several (say all) critical points bifurcate at $\lambda_0$:  is it then  possible to approximate $
\lambda_0$ by parameters such that the corresponding critical points are preperiodic? 
Of course in this question one has to discard a few ``trivial''  obstructions , e.g. when 
$\dim(\Lambda)$ is too small, so that there are not enough degrees of freedom to hope for an independent behavior of the critical points. Still after excluding these counterexamples, the answer to this problem is 
``no'' (see~\cite[Example 6.13]{prepercrit}),     the fundamental reason for this being  the failure of Montel's 
theorem in higher dimension. Using currents is a known way of circumventing this problem in higher dimensional 
dynamics, and,  as a matter of fact, the following theorem holds:

\begin{thm}[Bassanelli-Berteloot \cite{basber2}, Dujardin-Favre \cite{prepercrit}] \label{thm:inclusion}
Let $(f_\lambda)_{\lambda\in\Lambda}$ be a holomorphic family of rational maps of 
degree $d\geq 2$. Then for every $k\leq\dim(\Lambda)$, 
\begin{equation}\label{eq:inclusion}
\supp(\tbif^k)\subset \overline{ \set{\lambda,\ f_\lambda\text{ admits }k\text{ periodic critical points}}}. 
\end{equation}
\end{thm}

\noindent (This result was actually not stated explicitly  in~\cite{basber2, prepercrit},  see~\cite{survey} for this formulation. The converse inclusion is studied below.)

When $\Lambda$ is the moduli space $\poly_d$ of polynomials of degree $d$ 
with marked  critical points (which is a finite quotient of $\C^{d-1}$)  or the moduli space $\modu_d$ 
of rational maps of degree $d$ with marked critical points (which is   of dimension $2d-2$), 
we define the \emph{bifurcation measure} $\mubif$ to be the maximal 
exterior power of $\tbif$, that is  $\mubif = \tbif^{d-1}$ or $\mubif = \tbif^{2d-2}$, respectively. The following 
neat dynamical  characterization of $\supp(\mubif)$ can be obtained: 

\begin{thm}[Dujardin-Favre \cite{prepercrit}, Buff-Epstein\cite{buff-epstein}]\label{thm:pcf}
For $\Lambda = \poly_d$ or $\modu_d$, the support of 
$\mubif$ is the closure of (non-Lattès) strictly post-critically finite parameters, that is, 
parameters for which all critical points are preperiodic to a repelling cycle.
\end{thm}

 A version of this result for intermediate powers of $\tbif$ was   
  obtained in~\cite{higher}, which explains to what extent the converse inclusion in~\eqref{eq:inclusion} holds. 

\begin{proof}[Sketch of proof]
The most delicate point is to show that any non-Lattès post-critically finite parameter $\lambda_0$
  belongs to $\supp(\mubif)$.   To fix the ideas, assume that  $\Lambda = \mathcal M_d$. 
Observe that  $\lambda_0$ is an intersection point of a family of $(2d-2)$ hypersurfaces  of the form 
$$\set{\lambda\in \mathcal M_d, \ f^n_\lambda(c_i(\lambda))  = f^{n+k}_\lambda(c_i(\lambda))}$$ (one for each critical point). 
The proof in~\cite{buff-epstein} is based on two important ideas. 
The first one  consists in  proving that these hypersurfaces are smooth and transverse at $\lambda_0$: this is based  on
 Teichmüller-theoretic ideas. 
Then, using this transversality,    a version of Tan Lei's transfer principle between dynamical and parameter space
allows  to compare the mass of 
$\mubif$ in a carefully  scaled small polydisk about $\lambda_0$ with the mass of $\mu_{f_{\lambda_0}}$ near the $f^n(c_i)$, and conclude that this mass is positive.
\end{proof}  
  
In the space of polynomials of degree $d$, 
Theorem~\ref{thm:pcf}, together with   other   characterizations of $\supp(\mubif)$, 
e.g. in terms of landing of parameter rays, makes 
$\supp(\mubif)$ the natural analogue of the boundary of the Mandelbrot set for polynomials of 
 higher degree. This  motivates an investigation of  its topological and geometric properties. 
  First, it is a compact set, which,  for $d\geq 3$, 
   is strictly contained in the boundary of the   locus 
 $\mathcal C_d$ of polynomials with connected Julia set.  A topological consequence of Theorem~\ref{thm:inclusion} is that $\supp(\mubif)$ is contained in the closure of $\Int(\mathcal C_d)$; on the other hand it is unknown whether $\mathcal C_d$ is the closure of its interior. 
 Gauthier~\cite{gauthier_dimension} extended Shishikura's theorem to show that \emph{$\supp(\mubif)$ has maximal  Hausdorff dimension at each of its  points}. Let us also note that by using 
 advanced non-uniform hyperbolicity techniques, 
 it was shown   by Astorg, Gauthier, Mihalache and Vigny~\cite{AGMV} that  
  in the space $\modu_d$  of rational maps of degree $d$, 
 \emph{$\supp(\mubif)$ has positive volume}. 
 
The technical core of Theorems~\ref{thm:inclusion} and~\ref{thm:pcf} is the fact that $\tbif$ and its
 exterior powers
describe the asymptotic distribution of families  of dynamically defined hypersurfaces in parameter space, like 
parameters with a preperiodic critical point, or parameters with a periodic point of a given multiplier. Initiated in 
\cite{basber1, basber2, prepercrit}, this research theme has gradually evolved in  scope and sophistication, 
notably through its connections with arithmetic equidistribution (see~\cite{favre-gauthier_equidistribution}). 
 
 A striking and unexpected consequence of this technology is   
an asymptotic estimate for the number of hyperbolic components in $\modu_d$, which is so far not 
  accessible by other means. Recall that a \emph{hyperbolic component} is a connected component of the stability locus in which the dynamics is uniformly expanding on the Julia set. 
  We say that a hyperbolic component $\Omega$ is \emph{of disjoint type} $(n_1, \ldots, n_{2d-2})$ if the critical points are attracted by distinct attracting cycles of respective exact period $n_i$. 

\begin{thm}[Gauthier, Okuyama and Vigny \cite{GOV}]
The number $N(n)$ of hyperbolic components of disjoint type $(n, \ldots, n)$ in $\modu_d$ satisfies 
$$ N(n)\underset{n\to\infty}{\sim} \frac{d^{(2d-2)n}}{(2d-2)!} \int_{\modu_d} \mubif.$$
\end{thm}

\noindent  (An analogous formula holds for arbitrary disjoint 
type $(n_1, \ldots, n_{2d-2})$.) Note that the corresponding result in $\poly_d$ is much easier and follows essentially from Bézout's theorem 
(together with a transversality argument). The   value of $\int_{\modu_d} \mubif$ is  known only for $d=2$ \cite{GOV}. 

Once the bifurcation measure  is constructed on $\poly_d$ or $\modu_d$, it is natural to 
inquire about the dynamics of a $\mubif$-typical parameter. In $\modu_d$ this question is completely
 open so far. 
For the family of quadratic (and more generally unicritical) polynomials, it was shown by Graczyk-Swiatek~\cite{graczyk-swiatek} 
and  Smirnov~\cite{smirnov} in the late 90's that a $\mubif$-typical  parameter satisfies  
 the Collet-Eckmann condition; 
 in particular the local geometry of its Julia set is well understood.  
 These results are based on 
 combinatorial techniques and the landing of external and parameter rays, and the
  method  carries over for degree $d$ polynomials (see~\cite[Thm. 10]{prepercrit}). 
Interestingly, a completely new approach to the results of \cite{graczyk-swiatek, smirnov} was recently found,
 which applies to arbitrary families of rational maps. 

\begin{thm}[De Thélin, Gauthier and Vigny \cite{dethelin-gauthier-vigny}]\label{thm:dtgv}
Let $(f_\lambda)_{\lambda\in \Lambda}$   be an algebraic family of rational maps of degree $d$  with a marked critical point $c(\lambda)$. Let  $T_c$ be the bifurcation current associated to $c$ and $\norm{T_c}$ be the associated total variation measure.
 Then for $\norm{T_c}$-a.e. $\lambda$,    
\begin{equation}\label{eq:dtgv}
\liminf_{n\to\infty} \abs{Df_\lambda^n(c(\lambda))}\geq \unsur{2}\log d >0.
\end{equation}
\end{thm}

For the unicritical family $z^d+\lambda$, this statement  is precisely the typicality of the Collet-Eckmann  expansion  property. 

\begin{proof}[Sketch of proof]
This is  an application of the techniques of~\S\ref{subs:lamin_applications}. We may assume that $\Lambda$ is of dimension 1, so that $T_c$ is just a positive measure on $\Lambda$. Consider the sequence of iterated graphs $\Gamma_{f^n(c)}$, parameterized by $\gamma_n:\lambda\mapsto (\lambda, f^n_\lambda(c(\lambda)))$. 
Then, as explained above, $T_c= \pi_*(\hat T\wedge [\Gamma_c])$, where $\pi:\Lambda\times \pu\to \Lambda$ is the first projection and $\hat T$ is the natural $\hat f$-invariant current in $\Lambda\times \pu$. Using the 
$\hat f$-invariance of $\hat T$, we infer    that 
$$T_c = \pi_*\lrpar{ d^{-n} [\Gamma_{f^n(c)}]\wedge \hat T},  \text{ and conversely }  
(\gamma_n)_* (T_c) = d^{-n} [\Gamma_{f^n(c)}]\wedge \hat T.$$
Since the     $\Gamma_{f^n(c)}$ are algebraic curves of uniformly  bounded genus, by the results of \S\ref{subs:approximation}, the   part 
$\lrpar{d^{-n} [\Gamma_{f^n(c)}]}^r$ of these curves made of disks of size $r$
 has mass 
$1- O(r^2)$, and  since $\hat T$ has continuous potential, by Theorem~\ref{thm:isect}
the intersection $d^{-n} [\Gamma_{f^n(c)}]\wedge \hat T$ is carried by 
$\lrpar{d^{-n} [\Gamma_{f^n(c)}]}^r$, up to a small error $\eta(r)$. 
But to fill up a set of measure $1-\eta(r)$ of $d^{-n} [\Gamma_{f^n(c)}]\wedge \hat T$, at least $c(r) d^n$ disjoint such  disks are required, and pulling them back by $\gamma_n$ we get a set of $c(r) d^n$ disjoint 
disks in $\Lambda$, covering 
a set of measure  $1-\eta(r)$ for $T_c$, each of which mapped under $\gamma_n$ to a disk of size $r$. 
Being disjoint, most of the pulled-back disks in  $\Lambda$ have   area at most $Cd^{-n}$, 
so the derivative of $\gamma_n$ there must typically be larger than $Cd^{n/2}$. Analyzing how  the 
derivative of $\gamma_n$ is expressed in terms of the $Df_\lambda^k(c(\lambda))$, for $0\leq k\leq n$, 
 finally leads to~\eqref{eq:dtgv}.
\end{proof}

As already mentioned, the theory of bifurcation currents has deep connections with arithmetic dynamics, and 
related  rigidity problems in moduli spaces. A typical 
problem in this context is the classification of  families with a marked point  
$(f_\lambda, a(\lambda))$ for which the bifurcation current $T_a$ is ``abnormally regular''. 
The reader is referred to  the  
recent monograph \cite{favre-gauthier} by Favre and Gauthier for more on this topic.

\subsection{Stability/bifurcation theory in higher dimension} 

Moving to higher dimension, it is tempting to imitate the definition of $J$-stability 
by coining a definition of stability from the non-collision of periodic points.   
An obvious difficulty is that in this context the automatic  extension of holomorphic motions  fails and 
the relevance of this definition  needs to be justified, for instance by proving its equivalence with  
 other natural  ones. Due to the variety of possible situations, in higher dimension  the details depend on the category of maps under study. So far, this program has been fulfilled   in two cases: 
 polynomial automorphisms of $\cd$ (by Lyubich and the author), and holomorphic maps on $\pk$ (by Berteloot, Bianchi and Dupont). 

\subsubsection{Polynomial automorphisms of $\cd$} 
For a polynomial automorphism  $f$ of $\cd$ we can define  Julia sets $J^+$ and $J^-$ respectively 
associated to forward and backward iteration, as well as the ``small Julia set''     $J = J^+\cap J^-$, 
 and $J^*\subset J$ the closure of the set of 
saddle periodic points, which is also the support of the maximal entropy measure~\cite{bls}. Following~\cite{tangencies}, we say that a holomorphic family $(f_\lambda)_{\lambda\in \Lambda}$ 
of polynomial automorphisms of fixed dynamical degree $d$
  is \emph{weakly $J^*$-stable}   if (i) its  saddle points  do not bifurcate, hence (under mild assumptions) 
  so do all periodic points. (Here the numbering of  properties corresponds to that of the 1-dimensional case on page~\pageref{conditions}.)
  Then the holomorphic motion of saddle points extends to a \emph{branched holomorphic motion} of $J^*$ and 
  the condition is equivalent to (ii) $\lambda\mapsto J^*(f_\lambda)$ is continuous. Furthermore the branched 
  holomorphic motion extends to the ``big Julia set'' $J^+\cup J^-$. It remains an open question whether 
weak $J^*$-stability 
yields  a conjugacy on $J^*$ or $J$ (that is, whether an analogue of (iii) holds). 
It is proved in \cite{hyperbolic} that weak $J^*$-stability  implies a probabilistic form  of 
structural stability, that is, a conjugacy can be defined on a full measure 
subset for any hyperbolic measure. Also,   weak $J^*$-stability preserves uniform 
hyperbolicity \cite{hyperbolic, saddle}, so the familiar concept of hyperbolic component makes sense 
in this setting. 

Even if strictly speaking  polynomial automorphisms have no critical points, 
the main issue in  \cite{tangencies} is about condition (iv) (stability of critical points). 
Indeed, it a popular  analogue of a prerepelling critical point for a 2-dimensional diffeomorphism 
 is a heteroclinic tangency, so we are looking for a characterization of stability in terms of
  (absence of) tangencies. It is well-known that in dissipative 
  dynamics, homoclinic tangencies yield bifurcations from  saddles to sources, and the main point of \cite{tangencies} is to find a mechanism for the converse implication.  The key is the phenomenon of 
  \emph{semi-parabolic implosion}. 
  
  Before moving  on to this topic, let us point 
  out that so far there is no theory of bifurcation currents for automorphisms of $\cd$. 
  
  \begin{question}
  For polynomial automorphisms of $\cd$, is 
 stability characterized 
 by the harmonicity of the Lyapunov exponents of the maximal entropy measure?  
 In other words,  does an  analogue of  condition (v) above hold? 
  \end{question}

\subsubsection{Semi-parabolic implosion and tangencies}
Parabolic implosion refers to a set of phenomena, discovered by Douady and Lavaurs, 
 occurring when unfolding a periodic point with a rational 
indifferent multiplier. To be specific, consider a family  of the form $$f_\lambda(z) = (1+\lambda) z+ z^2+\hot$$ in a 
neighborhood of the origin, for small $\lambda$. For $\lambda = 0$, the fixed point 0 admits a basin of attraction 
$\mathcal B$. Now   If $\lambda$ approaches the origin tangentially to the imaginary axis, we can track precisely 
how the parabolic basin $\mathcal B$ ``implodes'' by ``passing through the eggbeater'' created between two slightly repelling fixed points $p_\lambda = 0$ and $q_\lambda\approx -\lambda$. More precisely, 
for well-chosen $\lambda_n$, $f^n_{\lambda_n}$ converges 
locally uniformly in $\mathcal B$ to a non-constant 
 \emph{Lavaurs map} $\psi:\mathcal B\to \C$, depending on $(\lambda_n)$. 
Of course for $\lambda_n\equiv 0$, $\psi = 0$: in this sense the limiting dynamics of $f_\lambda$ as $\lambda\to 0$ is richer than that of $f_0$. This gives rise to a wealth of dynamical phenomena at a such a parabolic bifurcation, like the discontinuity of the Julia set or the birth of 
hyperbolic set of large Hausdorff dimension, which are instrumental in  Shishikura's theorem that 
the boundary of the Mandelbrot set has dimension 2. 

Bedford, Smillie and Ueda \cite{bsu} extended this analysis 
to the unfolding of a semi-parabolic fixed point of 
multiplicity 2 in $\cd$, that is, of the form 
\begin{equation}\label{eq:multiplicity2}
f_\lambda(z,w) = ((1+\lambda)z+z^2 + \hot, b_\lambda w+\hot), \text{ with } \abs{b_0}<1.
\end{equation} In this dissipative situation, as before the Lavaurs map is a limit of iterates of the form $f^n_{\lambda_n}$,  
 its domain is the  attracting basin  $\mathcal B$ of the origin, but 
its values are contained  a  curve: the \emph{repelling petal} of the semi-parabolic point. 
For polynomial automorphisms, 
this leads to a precise description of 
 the  discontinuity of the  Julia sets $J$ and  $J^+$  at $\lambda = 0$. (See also Bianchi~\cite{bianchi_implosion} for 
 some results about the   implosion of general parabolic germs.)

If $(f_\lambda)$ is an arbitrary family of dissipative polynomial automorphisms, semi-parabolic
 bifurcations (of possibly  arbitrary multiplicity) occur densely in the bifurcation locus by definition. 
 A mechanism   producing homoclinic tangencies from  
 semi-parabolic implosion was designed in~\cite{tangencies}. Besides the analysis of Lavaurs maps (which is not as precise as in the multiplicity 2 
 case~\eqref{eq:multiplicity2}), this involves a construction of ``critical points'' in semi-parabolic basins, 
 which by definition are 
 tangencies between unstable manifolds (associated to some given saddle point) and 
  the foliation of the basin by
 strong stable manifolds. Surprisingly, this construction is based on   Wiman's classical  theorem on entire functions of slow growth, and requires a stronger dissipativity condition: $\abs{\mathrm{Jac}(f_\lambda)}< d^{-2}$ (\emph{substantially dissipative} regime). Altogether we obtain the following theorem, which confirms a classical conjecture of Palis in this setting:
 
 \begin{thm}[Dujardin and Lyubich \cite{tangencies}]
 In a substantially dissipative family of polynomial automorphisms of $\cd$, parameters with homoclinic tangencies are dense in the bifurcation locus. 
 \end{thm}

It is expected that this result holds  without the substantial dissipativity assumption. Also, it is an open question whether quadratic tangencies are always created in this process. A positive answer 
would yield an interesting link with the quadratic family, and 
add further evidence to the universality of the Mandelbrot set. 

\subsubsection{Holomorphic maps on $\pk$} The case of families of 
holomorphic maps on $\pk$ was studied by Berteloot, Bianchi and Dupont in \cite{bbd}. Here, as in the 
one-dimensional case, one starts with the stability of repelling periodic 
points. More precisely, one has to restrict to repelling 
points contained in the ``small  Julia set'' $J^*$ (which by definition is 
 the support of the maximal entropy mesure $\mu$), 
 since there can be a number of ``spurious'' repelling points outside 
 $J^*$. Then  Berteloot, Bianchi and Dupont obtain an almost complete generalization of the results of Ma\~né-Sad-Sullivan, Lyubich and DeMarco
 (that is, conditions (i) to (v) of pp. \pageref{conditions}-\pageref{condition5}). As before, a remaining issue is    whether this notion of weak 
$J^*$-stability implies structural stability on $J^*$. A main difference with the 1-dimensional case is that the 
characterization of bifurcation in terms of currents is now essential to establish the equivalence between the remaining conditions. More precisely the link between the instability of critical orbits and that of periodic points is 
provided by a formula à la Manning-Przytycki for the Lyapunov exponent of the maximal entropy measure. 
 
 We saw in  Theorems~\ref{thm:inclusion} and~\ref{thm:pcf} that the higher bifurcation currents $\tbif^k$ describe accurately certain higher codimensional phenomena in parameter space. It seems that the distinction between $\tbif$ and its powers is not as clear  in higher dimensional dynamics: in a recent work, Astorg and Bianchi \cite{astorg-bianchi} showed that in  a large portion of 
  the family of polynomial skew products of $\cd$,  the supports of all  currents $\tbif^k$ coincide with the bifurcation locus.  So the significance of  these higher bifurcation currents in this context is yet to be explored.

  \subsection{Robust bifurcations} As said before, due to the finiteness of the critical locus, 
  one-dimensional polynomial  and rational  maps are generically stable. 
Intuition from real dynamics suggests that this is not anymore the case in higher dimension.
As in the previous paragraph, we  discuss separately the cases of    polynomial automorphisms  and of 
holomorphic maps on $\P^k$. 

\subsubsection{Polynomial automorphisms} Given the characterization  of weak $J^*$-stability 
in~\cite{tangencies}, a straightforward adaptation of the one-dimensional argument for the density of stability 
shows that in any holomorphic family $(f_\lambda)$ of polynomial automorphisms of $\cd$, the union of (weakly $J^*$-)stable parameters together with 
 parameters with infinitely many sinks is dense. Prior to~\cite{tangencies}, it 
  was actually already known that stability is not a dense phenomenon in this context, due to the following remarkable result: 
  
  \begin{thm}[Buzzard \cite{buzzard}] 
  There exists $d>1$ and an open subset    
  $\Omega\subset \mathrm{Aut}_d(\cd)$  contained in the bifurcation locus. In particular maps with infinitely many sinks are dense in~$\Omega$.  
  \end{thm}
  
  Here $\mathrm{Aut}_d(\cd)$ is the space of polynomial automorphisms of $\cd$ of degree $d$. 
This deep theorem is nothing but the adaptation to the complex setting  of  
Newhouse's theorem (see~\cite{newhouse}) on the existence of surface 
diffeomorphisms with persistent homoclinic tangencies. It is obtained by first constructing transcendental examples and then 
approximating them  by polynomial ones, hence the degree $d$ is unknown and presumably very large.
The existence of this complex Newhouse phenomenon in arbitrary degree is a major open problem. 

\begin{question}
Is the bifurcation locus of non-empty interior in $\mathrm{Aut}_d(\cd)$ for any   $d\geq 2$?
\end{question}
 
As in the real case (cf. \cite{newhouse}), one may even expect that  robust bifurcations 
(that is, interior points of the bifurcation locus) are dense in the bifurcation locus, 
at least in the dissipative regime. 
For this, it is tempting to imitate  the  approach  of Shishikura's theorem on 
the Hausdorff dimension of $\fr M$
 and  use semi-parabolic implosion to 
 construct large bifurcation sets from a single parabolic bifurcation: in this sense the 
 density of robust bifurcations would be the  optimal generalization of Shishikura's theorem  
 to automorphisms of $\cd$. An interesting
  first  step would be to show that the bifurcation locus has maximal Hausdorff dimension at every point. More advanced techniques will certainly be needed to get open subsets: 
  an ambitious research program on 
  the intersection of complex Cantor sets was initiated by 
    Araujo, Moreira and Zamudio towards  this perspective (see~\cite{araujo-moreira, araujo-moreira-zamudio}).  

Biebler observed in~\cite{biebler_automorphisms} that the existence of robust bifurcations is actually more 
tractable in higher dimensions  and showed that: 
\emph{for every  $d\geq 2$ the bifurcation locus has non-empty interior in   
  $ \mathrm{Aut}_d(\C^3)$.} This is based on a distinct mechanism for robust bifurcation: the \emph{blenders}  of   
  Bonatti and Diaz \cite{bonatti-diaz}.  These are dynamically defined Cantor set which are   
  so fat in a certain ``direction''  that they intersect an open set of curves. The point of~\cite{biebler_automorphisms} is to use this feature as  
   a building block for persistent tangencies.

Finally, let us point out a recent beautiful result by    Yampolsky and Yang~\cite{yampolsky-yang}: \emph{the 
one-dimensional family of degree 2 Hénon maps
 with   a golden mean Siegel disk}
 $$f_a(x,y) = (x^2+c_a - ay, x), \text{ with } c_a= (1+a)\lrpar{\frac{\mu}{2}+\frac{a}{2\mu}}- \lrpar{\frac{\mu}{2}+\frac{a}{2\mu}}^2\text{ and }\mu   = e^{\pi(1+\sqrt{5})i}, $$
  \emph{is structurally unstable at every parameter with     small enough Jacobian $\abs{a}$}.
  This relies on  a completely different approach to persistent tangencies, based on Siegel renormalization.

\subsubsection{Holomorphic maps on $\P^k$} From the work of Berteloot, Bianchi and Dupont, 
we know that the basic phenomenon 
responsible for bifurcations  for holomorphic maps on $\pk$ 
is when the post-critical set intersects the small Julia set 
$J^*$. Thus, to obtain robust bifurcations, it is enough to  find a 
mechanism ensuring  a robust intersection between the post-critical set and $J^*$. 
A convenient tool for this is the Bonatti-Diaz blender, which leads to:

\begin{thm}[Dujardin \cite{robust}] For every $k\geq 2$ and $d\geq 2$, the bifurcation locus has non-empty interior in 
$\mathrm{Hol}_d(\P^k)$.  
\end{thm}

Here, $\mathrm{Hol}_d(\P^k)$ is the space of holomorphic maps on $\pk$ of degree $d$. 
A specific one-dimensional family of holomorphic maps of $\pd$ with a full bifurcation locus was found independently 
 by Bianchi and Taflin \cite{bianchi_taflin}. After this result, a natural question is that of the abundance of robust bifurcations in $\mathrm{Hol}_d(\P^k)$. Taflin \cite{taflin} showed that robust bifurcations are abundant near product polynomial maps of $\cd$, and Biebler 
 \cite{biebler_lattes} showed that Lattès maps of sufficiently large degree are accumulated by robust bifurcations. Blenders are involved directly or indirectly in both cases, and seem to appear quite naturally when a repelling   periodic point bifurcates to a saddle. Still, the general picture remains  elusive. 

 \begin{question}
 Is the bifurcation locus in $\mathrm{Hol}_d(\P^k)$ the closure of its interior? 
 \end{question}

Lastly, a
celebrated theorem   of McMullen asserts that any   stable algebraic families of rational maps on $\pu$ is either
 isotrivial or a family of flexible Lattès examples \cite{mcmullen_algorithms}. 
 Extending this result to higher dimensions is a promising  research problem; one main obstacle is that part of the argument relies on Thurston's topological characterization of rational functions. Related preliminary results have been obtained by Gauthier and Vigny \cite{gauthier-vigny_northcott}.

\section{(Non-)Wandering Fatou components}\label{sec:wandering}~

 The classification of Fatou components is a basic chapter of holomorphic dynamics. For rational maps 
in dimension 1, 
 periodic Fatou components can be classified into attracting basins, parabolic basins, and rotation 
domains (Siegel disks and Herman rings). The crowning achievement  of this  classification is the celebrated 
non-wandering domain theorem of Sullivan~\cite{sullivan}: \textit{for a one-dimensional rational map, any Fatou 
component is preperiodic}. 

In higher dimensions,  techniques from geometric function theory may be applied to classify periodic Fatou 
components. It is convenient to distinguish between recurrent and non-recurrent periodic 
components: a fixed Fatou 
component $\Omega$ is \emph{recurrent} if for some $x\in \Omega$, the $\omega$-limit set $\omega(x)$ is not 
completely contained in $\fr\Omega$. Recurrent Fatou components were classified 
in various classes of rational maps in~\cite{bs2, fornaess-sibony_classification,  Ueda_fatou, fornaess-rong}. The 
upshot is that in such a component either there is an transversely 
attracting submanifold (possibly  a point) or the dynamics is of rotation type. The situation is far less understood 
in the non-recurrent case. A notable exception is that of
 substantially  dissipative automorphisms of $\cd$, 
for which it was  shown by Lyubich and Peters~\cite{lyubich-peters}  that any non-recurrent Fatou component is 
the basin of a semi-parabolic periodic point. 

On the other hand it is immediately clear that the quasiconformal techniques used in Sullivan's proof 
are not generalizable to higher dimension. As it turns out, wandering components do 
exist in 2-dimensional polynomial dynamics:

\begin{thm}[Astorg, Buff, Dujardin, Peters and Raissy~\cite{ABDPR}]\label{thm:ABDPR}
If $0<a<1$ is sufficiently close to 1, the polynomial mapping of $\cd$ defined by:
\[ f: (z,w)  \longmapsto (p(z,w), q(w) )= \lrpar{z+z^2+a z^3 + \frac{\pi^2}{4}w, w-w^2} \]
admits a wandering Fatou component. 
\end{thm}

 The proof  is based on an original idea of M. Lyubich, and relies on a skew product version of parabolic implosion. It was further  implemented in other  situations  in~\cite{astorg-thaler-peters, hahn-peters}.

\begin{proof}[Sketch of proof] 
Write $p(z,w)   = p_0(z)+\e(z,w)$, with $p_0(z) = z+z^2$ and $\e(z,w)$ is thought of as a perturbative term.  
Start with an initial point $(z_0, w_0)$     such that $z_0$ belongs to the parabolic basin of attraction  of 
0 for $p_0$ and $w_0$ a small positive number, and let as usual 
$(z_n, w_n)  = f^n(z_0, w_0)$. Then $w_n = q^n(w)$ converges to 0 along the positive real axis, and 
$p^n_0(z_0)$ converges to 0 along the negative real axis. Therefore $z_n = p^n_0(z_0) +\e_n$ is pushed a little
faster towards the origin by the term $\e_n$. The terms in  $\e(z,w)$ are crafted  so 
that if $z_0$  is chosen carefully in some open set of initial conditions, 
the iterates $z_n$ indeed  pass the origin by going  ``through the eggbeater'' and come back close 
 to their initial position. So we can repeat this process and conclude that 
 $(z_0, w_0)$ belongs to some Fatou component. But since the returning time  increases with the number of 
 iterations, this Fatou component is not periodic, and we are done. 
 \end{proof}

At this stage the following natural questions arise:

\begin{question}~\label{question:wandering}
\begin{enumerate}
\item Are there other dynamical mechanisms leading to wandering Fatou components?
\item Find substantial families of higher dimensional rational mappings  without wandering domains.
\end{enumerate}
\end{question}

 Regarding   the first question, a 
mechanism for constructing wandering domains  in 2-dimensional smooth dynamics, based on the Newhouse 
phenomenon, was devised by Colli and Vargas~\cite{colli-vargas}. 
Berger and Biebler recently proved that this mechanism can be implemented in 
 certain 5-dimensional families of Hénon maps, leading to the following stunning  theorem:
 
 \begin{thm}[Berger and  Biebler~\cite{berger-biebler}]\label{thm:berger-biebler}
 There exists a polynomial automorphism of $\cd$ of degree 6 with a wandering Fatou component. 
 \end{thm}
 
This solves the existence problem for wandering Fatou components  for plane polynomial automorphisms, 
which does not seem to be amenable to the techniques of \cite{ABDPR}.

For    the second question, it is a classical fact that hyperbolic dynamics prevents the existence of 
wandering domains. Besides this observation, not much is known. In view of  
 Theorem~\ref{thm:ABDPR}, it is natural 
to  investigate  the case of skew products with a fixed attracting fiber, that is, of the form 
\begin{equation}\label{eq:skew}
 f(z,w) = (p(z), q(z,w)) \text{, with }p(0) = 0\text{ and }\abs{p'(0)}<1.
 \end{equation} 
 In this case it could be expected that 
Sullivan's theorem, together with the attracting nature of the invariant fiber should be enough 
to prevent the existence of wandering domains. 
 Embarrassingly enough, even in such a simple situation, there is no definitive answer so far, and furthermore 
 it was shown by Peters and Vivas~\cite{peters-vivas} that the above  naive intuition does not lead to a proof. 
 Here is the current status of the problem:
 
 \begin{thm}[Lilov, Peters-Smit, Ji]
 If $f$ is an attracting skew product  as in~\eqref{eq:skew}, then there are no wandering components near the 
 attracting fiber, whenever:
 \begin{itemize}
 \item   $p'(0)= 0$~\cite{lilov} or    more generally if $\abs{p'(0)}$ is small enough (with respect to $p$ and $q$)~\cite{ji1};
 \item $\abs{p'(0)}<1$ and $q(0, \cdot)$ satisfies some non-uniform hyperbolicity properties~\cite{peters-smit, ji2}. 
 \end{itemize} 
 \end{thm}

 There is currently no hope for a general understanding of the problem of wandering Fatou components in 
 several dimensions, and even going beyond skew products seems to be  a serious challenge. An interesting first 
 case to be considered is that of Fatou components in the  neighborhood of  an invariant  super-attracting line, 
 which would  cover for instance the case of regular polynomial mappings of $\cd$ near the line at infinity.


%

\bibliographystyle{abbrv} 
\bibliography{icm-bib}
 
\end{document}